\makeatletter

\documentclass{paper}
\usepackage[utf8]{inputenc}
\usepackage{amsmath}
\usepackage{amsthm}
\usepackage{amssymb}
\usepackage{graphicx}
\usepackage[numbers]{natbib}
\usepackage[unicode=true,
 bookmarks=false,
 breaklinks=false,pdfborder={0 0 1},backref=section,colorlinks=false]
 {hyperref}

\makeatletter
\newenvironment{lyxcode}
	{\par\begin{list}{}{
		\setlength{\rightmargin}{\leftmargin}
		\setlength{\listparindent}{0pt}
		\raggedright
		\setlength{\itemsep}{0pt}
		\setlength{\parsep}{0pt}
		\normalfont\ttfamily}%
	 \item[]}
	{\end{list}}

\usepackage{stmaryrd}
\makeatother

\begin{document}
\title{Modelisation of a rolling disk with Sympy}
\author{Luc Jaulin}
\institution{Lab-Sticc, ENSTA-Bretagne}

\maketitle
\textbf{Abstract}. This paper proposes a Lagrangian approach to find
the state equations of a disk rolling on a plane without friction.
The approach takes advantage of a symbolic computation to simplify
the reasoning. 

\section{Introduction}

\begin{figure}
\centering\includegraphics[width=9cm]{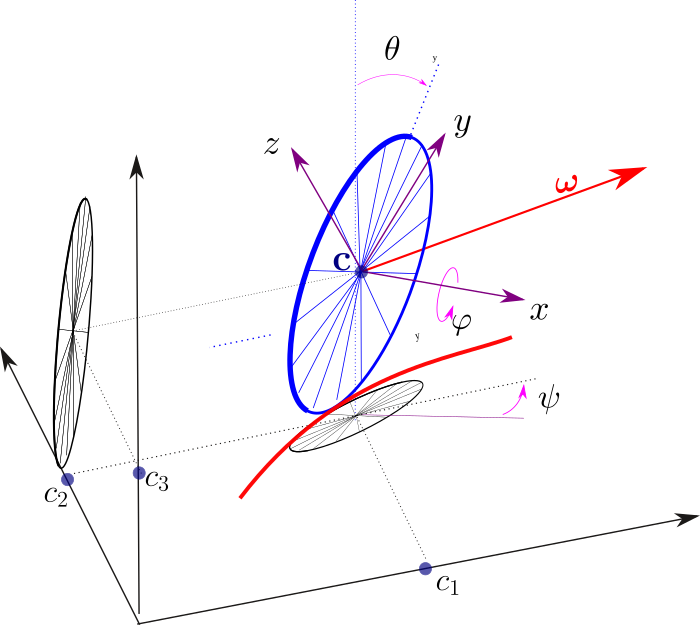}\caption{Disk (blue) rolling on a plane. The vertical and horizontal projections
are painted black}
\label{fig:rollingdisk}
\end{figure}

Consider a disk rolling on a plane without friction not sliding, as
shown on Figure \ref{fig:rollingdisk}.

We assume that the disk mass is $m=5\text{kg}$ and its radius is
$r=1m$. The gravity is taken as $g=9.81ms^{-2}.$

In this paper, we want to find the state equation describing the motion
of the disk. This problem has already been solved for over a century
(see \emph{e.g.} \citep{Appel:1900}) and even highly studied since
(see \emph{e.g.} \citep{OReilly:disk}). Extension to more general
wheel based vehicle have also been proposed (see \emph{e.g.} \citep{BoyerMauny18}
for the bicycle).

Computing the state equation for a rolling disk is a tedious task.
This paper takes advantage of symbolic computing (here the \texttt{sympy}
package of \texttt{Python}) in order to derive these state equations.
The Lagrangian approach \citep{wells67}, often applied to model robots
\citep{Corke11} will be chosen. 

\section{Modelisation}

\subsection{State vector}

We take the state vector $\mathbf{x}=(c_{1},c_{2},\varphi,\theta,\psi,\dot{\varphi},\dot{\theta},\dot{\psi})$
where $(c_{1},c_{2})$ is the vertical projection of center $\mathbf{c}=(c_{1},c_{2},c_{3})$
of the disk and $\varphi,\theta,\psi$ are the three Euler angle.
As illustrated by Figure \ref{fig:rollingdisk},
\begin{itemize}
\item $\varphi$ is the spin angle
\item $\theta$ is the stand angle, \emph{i.e.}, when $\theta=0$, the disk
is vertical
\item $\psi$ is the heading, \emph{i.e.}, the horizontal orientation of
the disk.
\end{itemize}
To use \texttt{sympy}, we first declare the symbolic variables and
functions to be used:
\begin{lyxcode}
t=symbols('t')~~~~~

m,g,r=~symbols('m~g~r')~~~~~

c1,c2~=~Function('c1')('t'),Function('c2')('t')~~~~~

dc1,dc2~=~Function('dc1')('t'),Function('dc2')('t')~~~~~

ddc1,ddc2~=~Function('ddc1')('t'),Function('ddc2')('t'),~~~~~

$\psi$,$\theta$,$\varphi$~=~Function('$\psi$')('t'),Function('$\theta$')('t'),Function('$\varphi$')('t')~~~~~

d$\psi$,d$\theta$,d$\varphi$~=~Function('d$\psi$')('t'),Function('d$\theta$')('t'),Function('d$\varphi$')('t')~~~~~

dd$\psi$,dd$\theta$,dd$\varphi$~=~Function('dd$\psi$')('t'),Function('dd$\theta$')('t'),Function('dd$\varphi$')('t')

$\varphi$,$\theta$,$\psi$~=~Function('$\varphi$')('t'),Function('$\theta$')('t'),Function('$\psi$')('t')~~~~~

d$\varphi$,d$\theta$,d$\psi$~=~Function('d$\varphi$')('t'),Function('d$\theta$')('t'),Function('d$\psi$')('t')~~~~~

dd$\varphi$,dd$\theta$,dd$\psi$~=~Function('dd$\varphi$')('t'),Function('dd$\theta$')('t'),Function('dd$\psi$')('t')~~~~~

$\lambda$1,$\lambda$2~=~Function('$\lambda$1')('t'),Function('$\lambda$2')('t')
\end{lyxcode}

\subsection{Orientation}

The orientation of the disk is fixed by the three Euler angles $\varphi,\theta,\psi$.
The corresponding orientation matrix is

\begin{equation}
\mathbf{R}=\left(\begin{array}{c}
\cos\theta\cos\psi\\
\cos\theta\sin\psi\\
-\sin\theta
\end{array}\right.\begin{array}{c}
-\cos\varphi\sin\psi+\sin\theta\cos\psi\sin\varphi\\
\cos\psi\cos\varphi+\sin\theta\sin\psi\sin\varphi\\
\cos\theta\sin\varphi
\end{array}\left.\begin{array}{c}
\sin\psi\sin\varphi+\sin\theta\cos\psi\cos\varphi\\
-\cos\psi\sin\varphi+\sin\theta\cos\varphi\sin\psi\\
\cos\theta\cos\varphi
\end{array}\right)
\end{equation}

It is built by the following \texttt{sympy} function
\begin{lyxcode}
def~Reuler($\varphi,\theta,\psi$):

~~~R$\varphi$~=~Matrix({[}{[}1,0,0{]},{[}0,cos($\varphi$),-sin($\varphi$){]},{[}0,sin($\varphi$),cos($\varphi$){]}{]})

~~~R$\theta$~=~Matrix({[}{[}cos($\theta$),0,sin($\theta$){]},{[}0,1,0{]},{[}-sin($\theta$),0,cos($\theta$){]}{]})~~~~~~~~~

~~~R$\psi$~=~Matrix({[}{[}cos($\psi$),-sin($\psi$),0{]},{[}sin($\psi$),cos($\psi$),0{]},{[}0,0,1{]}{]})~~~~~~~~~

~~~Return~~R$\psi${*}R$\theta${*}R$\varphi$
\end{lyxcode}
The rotation vector depends on the Euler angles and their derivatives.
Its expression \citep{jaulinISTEroben} can be obtained using the
relation 
\begin{equation}
\mathbf{R}^{\text{T}}\mathbf{\dot{R}}=\left(\begin{array}{ccc}
0 & -\omega_{r3} & \omega_{r2}\\
\omega_{r3} & 0 & -\omega_{r1}\\
-\omega_{r2} & \omega_{r1} & 0
\end{array}\right)
\end{equation}
which gives us the following \texttt{sympy} function
\begin{lyxcode}
def~wr(R):~

~~~W=Transpose(R){*}diff(R,t)~

~~~return~Matrix({[}{[}-W{[}1,2{]}{]},{[}W{[}0,2{]}{]},{[}-W{[}0,1{]}{]}{]})~~~~~~~~
\end{lyxcode}
We get
\begin{equation}
\boldsymbol{\omega}_{r}=\left(\begin{array}{c}
\dot{\varphi}-\dot{\psi}\sin\theta\\
\dot{\theta}\cos\varphi+\dot{\psi}\sin\varphi\cos\theta\\
-\dot{\theta}\sin\varphi+\dot{\psi}\cos\theta\cos\varphi
\end{array}\right).
\end{equation}

\subsection{Lagrangian}

In order, to get the state equation of the rolling disk, we use a
Lagrangian approach. For this, we need to express the Lagrangian $\mathcal{L}$
with respect to the state variables. Recall that 
\begin{equation}
\mathcal{L}=E_{K}-E_{p}
\end{equation}
where $E_{K}$ is the kinetic energy and $E_{p}$ is the potential
energy. We have

\begin{equation}
E_{K}=\frac{1}{2}\boldsymbol{\omega}_{r}^{\text{T}}\mathbf{I}\boldsymbol{\omega}_{r}+\frac{1}{2}m\|\dot{\mathbf{c}}\|^{2}
\end{equation}
where $\mathbf{c}=(c_{1},c_{2},c_{3})$ is the center of the disk
and 
\begin{equation}
\mathbf{I}=\left(\begin{array}{lll}
\frac{mr^{2}}{2} & 0 & 0\\
0 & \frac{mr^{2}}{4} & 0\\
0 & 0 & \frac{mr^{2}}{4}
\end{array}\right)
\end{equation}
is the inertia matrix of the disk. Moreover
\begin{equation}
E_{p}=mgr\cos\theta.
\end{equation}
Denote by
\begin{equation}
\mathbf{q}=\left(c_{1},c_{2},\varphi,\theta,\psi\right)
\end{equation}
the generalized coordinates of the system, \emph{i.e.}, the degrees
of freedom of the system. The Lagrangian, which is a function of $(\mathbf{q},\dot{\mathbf{q}})$,
is computed with \texttt{sympy} : 
\begin{lyxcode}
c1,c2,$\varphi$,$\theta$,$\psi$=list(q)~~~~~~~~~

c=Matrix({[}{[}c1{]},{[}c2{]},{[}r{*}cos($\theta$){]}{]})~~~~~~~~~

dc=diff(c,t)

R=Reuler($\varphi$,$\theta$,$\psi$)~~~~~~~~~

Ep=m{*}g{*}c{[}2{]}~~~~~~~~

I=Matrix({[}{[}1/2{*}m{*}r{*}{*}2,0,0{]},{[}0,1/4{*}m{*}r{*}{*}2,0{]},{[}0,0,1/4{*}m{*}r{*}{*}2{]}{]})

Ek=1/2{*}m{*}(dc1{*}{*}2+dc2{*}{*}2+dc3{*}{*}2)+(1/2){*}wr(R).dot(I{*}wr(R))~~~~~~~~~

L=Ek-Ep
\end{lyxcode}
which yields
\begin{equation}
\begin{array}{ccl}
\mathcal{L}(\mathbf{q},\dot{\mathbf{q}}) & = & \frac{1}{2}m\left(\dot{c}_{1}^{2}+\dot{c}_{2}^{2}+r^{2}\sin^{2}\theta\dot{\theta}^{2}\right)\\
 &  & +\frac{1}{8}mr^{2}\left(2\left(\dot{\varphi}-\sin\theta\dot{\psi}\right)^{2}+\dot{\theta}^{2}+\cos^{2}\theta\dot{\psi}^{2}\right)\\
 &  & -mgr\cos\theta
\end{array}
\end{equation}

The evolution of $\mathbf{q}$ obeys to the Lagrange's equation for
holonomic systems
\begin{equation}
\underset{\mathcal{Q}(\mathbf{q},\dot{\mathbf{q}},\ddot{\mathbf{q}})}{\underbrace{\frac{d}{dt}\left(\frac{\partial\mathcal{L}}{\partial\dot{\mathbf{q}}}\right)-\frac{\partial\mathcal{L}}{\partial\mathbf{q}}}}=\boldsymbol{\tau}\label{eq:EulerLagEq}
\end{equation}
where $\boldsymbol{\tau}$ are the constraint forces. The $i$th component
$\tau_{i}$ of $\boldsymbol{\tau}$ is associated to the $i$th component
$q_{i}$ of $\mathbf{q}$. Now
\begin{equation}
\frac{\partial\mathcal{L}}{\partial\dot{\mathbf{q}}}=\left(\begin{array}{c}
m\dot{c}_{1}\\
m\dot{c}_{2}\\
\frac{1}{2}mr^{2}(\dot{\varphi}-\dot{\psi}\sin\theta)\\
mr^{2}(\sin^{2}\theta+\frac{1}{4})\dot{\theta}\\
\frac{1}{4}mr^{2}(-2\dot{\varphi}\sin\theta+\dot{\psi}\sin^{2}\theta+\dot{\psi})
\end{array}\right)
\end{equation}
Thus
\begin{equation}
\frac{d}{dt}\left(\frac{\partial\mathcal{L}}{\partial\dot{\mathbf{q}}}\right)=\left(\begin{array}{c}
m\ddot{c}_{1}\\
m\ddot{c}_{2}\\
mr^{2}(\ddot{\varphi}-\ddot{\psi}\sin\theta-\dot{\theta}\dot{\psi}\cos\theta)/2\\
\frac{1}{4}mr^{2}(4\ddot{\theta}\sin^{2}\theta+\ddot{\theta}+4\dot{\theta}^{2}\sin(2\theta))\\
\frac{1}{4}mr^{2}(-2\ddot{\varphi}\sin\theta+\ddot{\psi}\sin^{2}\theta+\ddot{\psi}-2\dot{\theta}\dot{\varphi}\cos\theta+\dot{\theta}\dot{\psi}\sin(2\theta))
\end{array}\right)
\end{equation}
Moreover
\begin{equation}
\frac{\partial\mathcal{L}}{\partial\mathbf{q}}=\left(\begin{array}{c}
0\\
0\\
0\\
\frac{1}{4}mr\left(4g\sin\theta+2r\dot{\theta}^{2}\sin(2\theta)-2r\dot{\varphi}\dot{\psi}\cos\theta+\frac{r\dot{\psi}^{2}\sin(2\theta)}{2}\right)\\
0
\end{array}\right)
\end{equation}
The left hand side of the Euler Lagrange equation (\ref{eq:EulerLagEq})
is thus

\begin{equation}
\mathcal{Q}(\mathbf{q},\dot{\mathbf{q}},\ddot{\mathbf{q}})=\left(\begin{array}{c}
m\ddot{c}_{1}\\
m\ddot{c}_{2}\\
\frac{1}{2}mr^{2}(\ddot{\varphi}-\ddot{\psi}\sin\theta-\dot{\theta}\dot{\psi}\cos\theta)\\
\frac{1}{8}mr\left(-8g\sin\theta+8r\ddot{\theta}\sin^{2}\theta+2r\ddot{\theta}+4r\dot{\theta}^{2}\sin(2\theta)+4r\dot{\varphi}\dot{\psi}\cos\theta-r\dot{\psi}^{2}\sin(2\theta)\right)\\
\frac{1}{4}mr^{2}\left(-2\ddot{\varphi}\sin\theta+\ddot{\psi}\sin^{2}\theta+\ddot{\psi}-2\dot{\theta}\dot{\varphi}\cos\theta+\dot{\theta}\dot{\psi}\sin(2\theta)\right)
\end{array}\right)\label{eq:Q:q:dq:ddq}
\end{equation}

This expression for $\mathcal{Q}(\mathbf{q},\dot{\mathbf{q}},\ddot{\mathbf{q}})$
is obtained directly using \texttt{sympy} by 
\begin{lyxcode}
q=Matrix({[}c1,c2,$\varphi$,$\theta$,$\psi${]})~~~~~

dq=Matrix({[}dc1,dc2,d$\varphi$,d$\theta$,d$\psi${]})~~~~~

ddq=Matrix({[}ddc1,ddc2,dd$\varphi$,dd$\theta$,dd$\psi${]})

Q=diff(L.jacobian(dq),t)-L.jacobian(q)
\end{lyxcode}

\subsection{Non holonomic constraints}

\begin{figure}[h]
\centering\includegraphics[width=12cm]{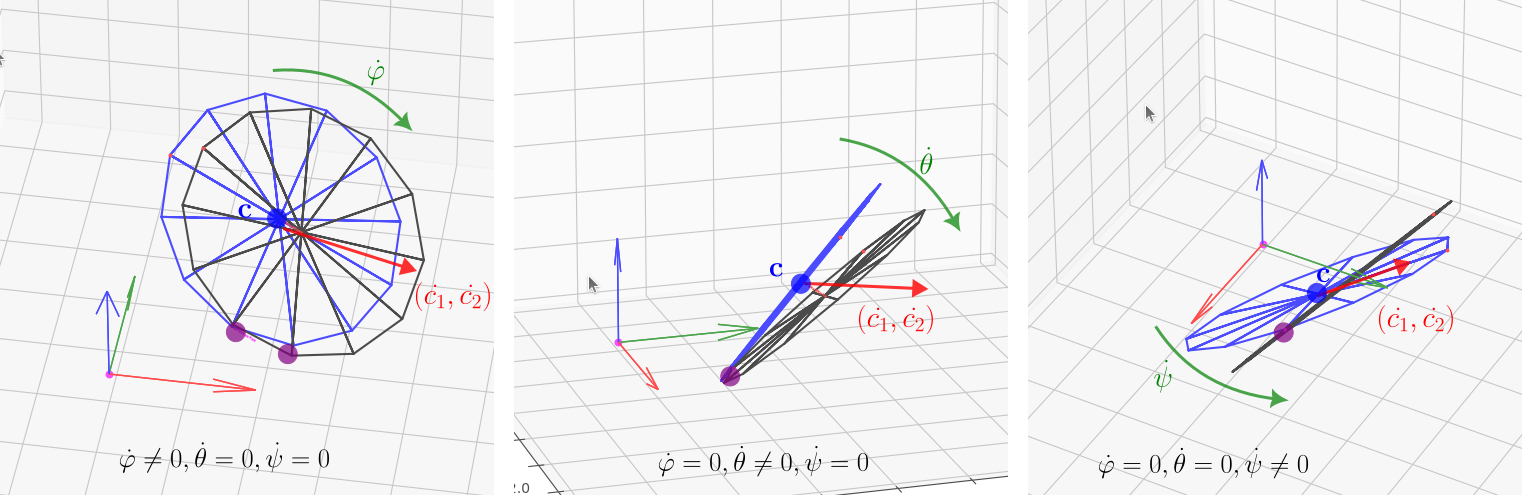}\caption{Deriving the non holonomic constraints}

\label{fig:dalembert.png}
\end{figure}

If the ground is a flat frozen lake where the disk can slide without
any friction in both direction (horizontally and laterally), then
the state vector is 
\begin{equation}
(\mathbf{q},\dot{\mathbf{q}})=(c_{1},c_{2},\varphi,\theta,\psi,\dot{c}_{1},\dot{c}_{2},\dot{\varphi},\dot{\theta},\dot{\psi}),
\end{equation}
\emph{i.e.}, the state is composed of the degrees of freedom $\mathbf{q}$
and they derivatives $\dot{\mathbf{q}}$. Now, in our case, no sliding
is possible and the disk can only roll. Due to this rolling constraint,
$\text{(\ensuremath{\mathbf{q}},\ensuremath{\dot{\mathbf{q}}})}$
are linked by some differential constraints. These constraints are
needed to derive the state equations with $\mathbf{x}=(c_{1},c_{2},\varphi,\theta,\psi,\dot{\varphi},\dot{\theta},\dot{\psi})$
as a state vector. Since we have two variables to eliminate (here
$\dot{c}_{1}$ and $\dot{c}_{2}$), we need to find two more differential
constraints which are generated by the ground forces. These constraints
translate the fact that the point of the disk in contact with the
ground has a zero velocity. It means that the disk neither slides
tangentially (first equation) nor laterally (second equation). We
understand that these two equations have the form
\begin{equation}
\begin{array}{ccc}
\dot{c}_{1} & = & \alpha_{1}\cdot\dot{\varphi}+\alpha_{2}\cdot\dot{\theta}+\alpha_{3}\cdot\dot{\psi}\\
\dot{c}_{2} & = & \beta_{1}\cdot\dot{\varphi}+\beta_{2}\cdot\dot{\theta}+\beta_{3}\cdot\dot{\psi}
\end{array}
\end{equation}
 where the $\alpha_{i}$'s and the $\beta_{i}'s$ depend on $\mathbf{q}$.
More precisely, this corresponds to the \emph{non holonomic} constraints
given by
\begin{equation}
\begin{array}{ccc}
\dot{c}_{1} & = & r\sin\psi\cdot\dot{\varphi}+r\cos\psi\cos\theta\cdot\dot{\theta}-r\sin\psi\sin\theta\cdot\dot{\psi}\\
\dot{c}_{2} & = & -r\cos\psi\cdot\dot{\varphi}+r\sin\psi\cos\theta\cdot\dot{\theta}+r\cos\psi\sin\theta\cdot\dot{\psi}
\end{array}\label{eq:nonholonomous:constraints}
\end{equation}
Figure \ref{fig:dalembert.png} illustrates how this formula is obtained.
The left subfigure shown that, when $\dot{\theta}=0$ and $\dot{\psi}=0$,
we have 
\begin{equation}
\begin{array}{ccc}
\dot{c}_{1} & = & r\sin\psi\cdot\dot{\varphi}\\
\dot{c}_{2} & = & -r\cos\psi\cdot\dot{\varphi}
\end{array}
\end{equation}
The subfigure in the center illustrates that if $\dot{\varphi}=0$
and $\dot{\psi}=0$, 
\begin{equation}
\begin{array}{ccc}
\dot{c}_{1} & = & r\cos\psi\cos\theta\cdot\dot{\theta}\\
\dot{c}_{2} & = & r\sin\psi\cos\theta\cdot\dot{\theta}
\end{array}
\end{equation}

The right subfigure illustrates that if $\dot{\theta}=0$ and $\dot{\varphi}=0$,
\begin{equation}
\begin{array}{ccc}
\dot{c}_{1} & = & -r\sin\psi\sin\theta\cdot\dot{\psi}\\
\dot{c}_{2} & = & r\cos\psi\sin\theta\cdot\dot{\psi}
\end{array}
\end{equation}
By superposition, we get Equation (\ref{eq:nonholonomous:constraints}).
These constraints are said to be non holonomic since they will not
allow us to express our system with a state composed of some degrees
of freedom $\mathbf{q}$ and their derivatives $\mathbf{\dot{q}}$. 

\subsection{D'Alembert's principle}

We need to find an expression for $\boldsymbol{\tau}$ which occurs
in the right hand side of the Euler-Lagrange equation (\ref{eq:EulerLagEq}).
The components for $\boldsymbol{\tau}$ correspond to the generalized
forces applied to our system. Now, in our specific case, we have no
friction no external forces to thrust or slow down the disk. The only
forces that apply are the reaction of the ground on the disk. 

In order to use this information, let us to recall the principle of
d'Alembert: \emph{for arbitrary virtual displacements, the constraint
forces don't do any work}. 

The virtual displacements are infinitesimal changes $\delta\mathbf{q}=\mathbf{q}(t+dt)-\mathbf{q}(t)$
for $\mathbf{q}(t)$ that should be consistent with some feasible
trajectories. For our rolling disk, the virtual displacements satisfy
\begin{equation}
\begin{array}{ccccc}
\delta c_{1}-r\sin\psi\cdot\delta\varphi-r\cos\psi\cos\theta\cdot\delta\theta+r\sin\psi\sin\theta\cdot\delta\psi & = & 0 & \, & (i)\\
\delta c_{2}+r\cos\psi\cdot\delta\varphi-r\sin\psi\cos\theta\cdot\delta\theta-r\cos\psi\sin\theta\cdot\delta\psi & = & 0 &  & (ii)
\end{array}\label{eq:virtual:disp}
\end{equation}
for the same reasons than those used to derive (\ref{eq:nonholonomous:constraints}).
The fact that there is no work translates into
\begin{equation}
\begin{array}{ccc}
\delta W=\boldsymbol{\tau}^{\text{T}}\delta\mathbf{q}=\tau_{c_{1}}\cdot\delta c_{1}+\tau_{c_{2}}\cdot\delta c_{2}+\tau_{\varphi}\cdot\delta\varphi+\tau_{\theta}\cdot\delta\theta+\tau_{\psi}\cdot\delta\psi & = & 0.\end{array}\label{eq:no:work}
\end{equation}

Equivalently, Equation (\ref{eq:no:work}) is a linear combination
of the two equations (\ref{eq:virtual:disp}), \emph{i.e.}, (\ref{eq:no:work})
$=\lambda_{1}\cdot$(\ref{eq:virtual:disp},i)$+\lambda_{2}\cdot$(\ref{eq:virtual:disp},ii)
and the $\lambda_{i}$ are the Lagrange parameters. Therefore:
\begin{equation}
\left(\begin{array}{c}
\tau_{c_{1}}\\
\tau_{c_{2}}\\
\tau_{\varphi}\\
\tau_{\theta}\\
\tau_{\psi}
\end{array}\right)=\underset{\boldsymbol{\tau}(\mathbf{q},\text{\ensuremath{\boldsymbol{\lambda}}})}{\underbrace{\lambda_{1}\cdot\left(\begin{array}{c}
1\\
0\\
-r\sin\psi\\
-r\cos\psi\cos\theta\\
r\sin\psi\sin\theta
\end{array}\right)+\lambda_{2}\cdot\left(\begin{array}{c}
0\\
1\\
r\cos\psi\\
-r\sin\psi\cos\theta\\
-r\cos\psi\sin\theta
\end{array}\right)}}\label{eq:tau:lam}
\end{equation}

Using (\ref{eq:Q:q:dq:ddq}) and (\ref{eq:tau:lam}), we get that
the Euler Lagrange equation rewrites into
\begin{equation}
\mathcal{Q}(\mathbf{q},\dot{\mathbf{q}},\ddot{\mathbf{q}})-\boldsymbol{\tau}(\mathbf{q},\text{\ensuremath{\boldsymbol{\lambda}}})=\mathbf{0}
\end{equation}
i.e.,
\begin{equation}
\begin{array}{ccc}
m\ddot{c}_{1}-\lambda_{1} & = & 0\\
m\ddot{c}_{2}-\lambda_{2} & = & 0\\
\frac{1}{2}mr^{2}(\ddot{\varphi}-\ddot{\psi}\sin\theta-\dot{\theta}\dot{\psi}\cos\theta)+\lambda_{1}r\sin\psi-\lambda_{2}r\cos\psi & = & 0\\
\frac{1}{8}mr\left(-8g\sin\theta+8r\ddot{\theta}\sin^{2}\theta+2r\ddot{\theta}+4r\dot{\theta}^{2}\sin(2\theta)+4r\dot{\varphi}\dot{\psi}\cos\theta-r\dot{\psi}^{2}\sin(2\theta)\right)\\
+\lambda_{1}r\cos\psi\cos\theta+\lambda_{2}r\sin\psi\cos\theta & = & 0\\
\frac{1}{4}mr^{2}\left(-2\ddot{\varphi}\sin\theta+\ddot{\psi}\sin^{2}\theta+\ddot{\psi}-2\dot{\theta}\dot{\varphi}\cos\theta+\dot{\theta}\dot{\psi}\sin(2\theta)\right)-\lambda_{1}r\sin\psi\sin\theta+\lambda_{2}r\cos\psi\sin\theta & = & 0
\end{array}\label{eq:eulerlag:eq}
\end{equation}

This system is made of 5 equations which are linear in 7 variables
: $\boldsymbol{\lambda}=(\lambda_{1},\lambda_{2})$ and $\ddot{\mathbf{q}}=(\ddot{c}_{1},\ddot{c}_{2},\ddot{\varphi},\ddot{\theta},\ddot{\psi}).$

In order to square the system, we need to add two equations (to get
7 equations). They can be derived from non-holonomic constraints (\ref{eq:nonholonomous:constraints})
given by 
\[
\underset{\mathbf{a}(\mathbf{q},\dot{\mathbf{q}})}{\underbrace{\underset{\mathbf{A}(\mathbf{q})}{\underbrace{\left(\begin{array}{ccccc}
1 & 0 & \,\,-r\sin\psi\,\, & \,\,-r\cos\psi\cos\theta\,\, & r\sin\psi\sin\theta\\
0 & 1 & r\cos\psi & -r\sin\psi\cos\theta & -r\cos\psi\sin\theta
\end{array}\right)}}\cdot\dot{\mathbf{q}}}}=\left(\begin{array}{c}
0\\
0
\end{array}\right)
\]
Let us differentiate this equation. We get:
\begin{equation}
\underset{=\frac{d}{dt}\mathbf{a}(\mathbf{q},\dot{\mathbf{q}})}{\underbrace{\frac{\partial\mathbf{a}(\mathbf{q},\dot{\mathbf{q}})}{\partial\mathbf{q}}\cdot\dot{\mathbf{q}}+\frac{\partial\mathbf{a}(\mathbf{q},\dot{\mathbf{q}})}{\partial\dot{\mathbf{q}}}\cdot\ddot{\mathbf{q}}}}=\mathbf{0}\label{eq:d_dta}
\end{equation}

We add these two equations to (\ref{eq:eulerlag:eq}) to get
\begin{equation}
\underset{=\mathcal{S}(\mathbf{q},\dot{\mathbf{q}},\ensuremath{\boldsymbol{\lambda}},\ddot{\mathbf{q}})}{\underbrace{\left(\begin{array}{c}
\frac{d}{dt}\mathbf{a}(\mathbf{q},\dot{\mathbf{q}})\\
\mathcal{Q}(\mathbf{q},\dot{\mathbf{q}},\ddot{\mathbf{q}})-\boldsymbol{\tau}(\mathbf{q},\text{\ensuremath{\boldsymbol{\lambda}}})
\end{array}\right)}}=\left(\begin{array}{c}
\mathbf{0}\\
\mathbf{0}
\end{array}\right)\label{eq:S}
\end{equation}

Now, $\frac{d}{dt}\mathbf{a}(\mathbf{q},\dot{\mathbf{q}})$ is linear
in $\ddot{\mathbf{q}}$ (see \ref{eq:d_dta}), $\mathcal{Q}(\mathbf{q},\dot{\mathbf{q}},\ddot{\mathbf{q}})$
is linear in $\ddot{\mathbf{q}}$, $\boldsymbol{\tau}(\mathbf{q},\text{\ensuremath{\boldsymbol{\lambda}}})$
is linear in $\text{\ensuremath{\boldsymbol{\lambda}}}$. Therefore
$\mathcal{S}(\ensuremath{\boldsymbol{\lambda}},\mathbf{q},\dot{\mathbf{q}},\ddot{\mathbf{q}})$
is linear in $\ensuremath{(\boldsymbol{\lambda}},\ddot{\mathbf{q})}$,
\emph{i.e.}, 
\begin{equation}
\mathcal{S}(\mathbf{q},\dot{\mathbf{q}},\ensuremath{\boldsymbol{\lambda}},\ddot{\mathbf{q}})=\mathbf{M}(\mathbf{q},\dot{\mathbf{q}})\left(\begin{array}{c}
\boldsymbol{\lambda}\\
\ddot{\mathbf{q}}
\end{array}\right)-\mathbf{b}(\mathbf{q},\dot{\mathbf{q}}).\label{eq:S:Mb}
\end{equation}
where $\mathbf{M}(\mathbf{q},\dot{\mathbf{q}})$ is called the \emph{mass
matrix}. Since it does not depend on only depends on $\dot{\mathbf{q}}$,
we will write $\mathbf{M}(\mathbf{q})$ instead of $\mathbf{M}(\mathbf{q},\dot{\mathbf{q}})$.
An expression for $\mathbf{M}(\mathbf{q})$ and $\mathbf{b}(\mathbf{q},\dot{\mathbf{q}})$
can be obtained from an expression for $\mathcal{S}$ by (see \ref{eq:S:Mb}):
\[
\begin{array}{ccc}
\mathbf{M}(\mathbf{q}) & = & \frac{\partial\mathcal{S}}{\partial(\ensuremath{\boldsymbol{\lambda}},\ddot{\mathbf{q}})}\\
\mathbf{b}(\mathbf{q},\dot{\mathbf{q}}) & = & \mathbf{M}(\mathbf{q})\left(\begin{array}{c}
\boldsymbol{\lambda}\\
\ddot{\mathbf{q}}
\end{array}\right)-\mathcal{S}(\mathbf{q},\dot{\mathbf{q}},\ensuremath{\boldsymbol{\lambda}},\ddot{\mathbf{q}})
\end{array}
\]

Thus, (\ref{eq:S:Mb}) rewrites into
\begin{equation}
\mathbf{M}(\mathbf{q})\left(\begin{array}{c}
\boldsymbol{\lambda}\\
\ddot{\mathbf{q}}
\end{array}\right)=\mathbf{b}(\mathbf{q},\dot{\mathbf{q}}),
\end{equation}
An expression for $\mathbf{M}(\mathbf{q})$ and $\mathbf{b}(\mathbf{q},\dot{\mathbf{q}})$
are obtained by
\begin{lyxcode}
A=Matrix({[}{[}1,0,-r{*}sin($\psi$),-r{*}cos($\psi$){*}cos($\theta$),r{*}sin($\psi$){*}sin($\theta$){]},

~~~~~~~~~~{[}0,1,~r{*}cos($\psi$),-r{*}sin($\psi$){*}cos($\theta$),-r{*}cos($\psi$){*}sin($\theta$){]}{]})~~~~~

$\tau$=$\lambda$1{*}A{[}0,:{]}+$\lambda$2{*}A{[}1,:{]}~

a=A{*}dq~~~~~

da=diff(a,t)~~~~~

S=Matrix({[}da,{*}list(Q-$\tau$){]})~~~~~

M=S.jacobian({[}$\lambda$1,$\lambda$2,ddq{]})~~~~~

b=M{*}Matrix({[}$\lambda$1,$\lambda$2,ddq{]})-S~~~~~

\end{lyxcode}
We get
\begin{equation}
\mathbf{M}=\left(\begin{array}{ccccccc}
0 & 0 & 1 & 0 & -r\sin\psi & -r\cos\theta\cos\psi & r\sin\theta\sin\psi\\
0 & 0 & 0 & 1 & r\cos\psi & -r\sin\psi\cos\theta & -r\sin\theta\cos\psi\\
-1 & 0 & m & 0 & 0 & 0 & 0\\
0 & -1 & 0 & m & 0 & 0 & 0\\
r\sin\psi & -r\cos\psi & 0 & 0 & \frac{mr^{2}}{2} & 0 & -\frac{mr^{2}\sin\theta}{2}\\
r\cos\theta\cos\psi & \,r\sin\psi\cos\theta & 0 & 0 & 0 & mr^{2}(\sin^{2}\theta+\frac{1}{4}) & 0\\
-r\sin\theta\sin\psi\, & r\sin\theta\cos\psi & 0 & 0 & \,-\frac{mr^{2}}{2}\sin\theta & 0 & \frac{mr^{2}(\sin^{2}\theta+1)}{4}
\end{array}\right)
\end{equation}
and 
\begin{equation}
\mathbf{b}=\left(\begin{array}{c}
r(-\dot{\theta}^{2}\sin\theta\cos\psi-2\dot{\theta}\dot{\psi}\sin\psi\cos\theta+\dot{\varphi}\dot{\psi}\cos\psi-\dot{\psi}^{2}\sin\theta\cos\psi)\\
r(-\dot{\theta}^{2}\sin\theta\sin\psi+2\dot{\theta}\dot{\psi}\cos\theta\cos\psi+\dot{\varphi}\dot{\psi}\sin\psi-\dot{\psi}^{2}\sin\theta\sin\psi)\\
0\\
0\\
mr^{2}\dot{\theta}\dot{\psi}\cos\theta/2\\
\frac{mr}{8}(8g\sin\theta-4r\dot{\theta}^{2}\sin(2\theta)-4r\dot{\varphi}\dot{\psi}\cos\theta+r\dot{\psi}^{2}\sin(2\theta))\\
\frac{mr^{2}}{2}(\dot{\varphi}-\dot{\psi}\sin\theta)\dot{\theta}\cos\theta
\end{array}\right)
\end{equation}
We isolate $\boldsymbol{\lambda}$ and $\ddot{\mathbf{q}}$ by 
\begin{lyxcode}
$\lambda$1,$\lambda$2,ddc1,ddc2,dd$\varphi$,dd$\theta$,dd$\psi$=list((M.inv(){*}b))

~
\end{lyxcode}
and we get

\begin{equation}
\begin{array}{ccc}
\lambda_{1} & = & m\frac{6g\sin(2\theta)\cos\psi-15r\dot{\theta}^{2}\sin\theta\cos\psi-5r\dot{\theta}\dot{\psi}\sin\psi\cos\theta+18r\dot{\varphi}\dot{\psi}\sin^{2}\theta\cos\psi-3r\dot{\varphi}\dot{\psi}\cos\psi-15r\dot{\psi}^{2}\sin^{3}\theta\cos\psi}{15}\\
\lambda_{2} & = & m\frac{6g\sin(2\theta)\sin\psi-15r\dot{\theta}^{2}\sin\theta\sin\psi+5r\dot{\theta}\dot{\psi}\cos\theta\cos\psi+18r\dot{\varphi}\dot{\psi}\sin^{2}\theta\sin\psi-3r\dot{\varphi}\dot{\psi}\sin\psi-15r\dot{\psi}^{2}\sin^{3}\theta\sin\psi}{15}\\
\ddot{c}_{1} & = & \frac{2g\sin(2\theta)\cos\psi}{5}-r\dot{\theta}^{2}\sin\theta\cos\psi-\frac{r}{3}\dot{\theta}\dot{\psi}\sin\psi\cos\theta+\frac{6}{5}r\dot{\varphi}\dot{\psi}\sin^{2}\theta\cos\psi-\frac{r}{5}\dot{\varphi}\dot{\psi}\cos\psi-r\dot{\psi}^{2}\sin^{3}\theta\cos\psi\\
\ddot{c}_{2} & = & \frac{2g\sin(2\theta)\sin\psi}{5}-r\dot{\theta}^{2}\sin\theta\sin\psi+\frac{r}{3}\dot{\theta}\dot{\psi}\cos\theta\cos\psi+\frac{6}{5}r\dot{\varphi}\dot{\psi}\sin^{2}\theta\sin\psi-\frac{r}{5}\dot{\varphi}\dot{\psi}\sin\psi-r\dot{\psi}^{2}\sin^{3}\theta\sin\psi\\
\ddot{\varphi} & = & 2\dot{\varphi}\dot{\theta}\tan\theta+\frac{5}{3}\dot{\theta}\dot{\psi}\cos\theta\\
\ddot{\theta} & = & \frac{4}{5r}g\sin\theta-\frac{6}{5}\dot{\varphi}\dot{\psi}\cos\theta+\frac{1}{2}\dot{\psi}^{2}\sin(2\theta)\\
\ddot{\psi} & = & \frac{2\dot{\varphi}\dot{\theta}}{\cos\theta}
\end{array}\label{eq:lam:ddq}
\end{equation}

Finally, the \texttt{Python} function associated to the evolution
equation of the rolling disk is
\begin{lyxcode}
lambdify((c1,c2,$\varphi$,$\theta$,$\psi$,d$\varphi$,d$\theta$,d$\psi$,m,g,r),

~~~~~~~(dc1-a{[}0{]},dc2-a{[}1{]},d$\varphi$,d$\theta$,d$\psi$,dd$\varphi$,dd$\theta$,dd$\psi$))

\end{lyxcode}
or equivalently, the state equation are
\begin{equation}
\begin{array}{ccccc}
\left(\begin{array}{c}
\dot{c}_{1}\\
\dot{c}_{2}
\end{array}\right) & = & r\left(\begin{array}{ccccc}
\sin\psi & \, & \cos\psi\cos\theta & \, & -\sin\psi\sin\theta\\
-\cos\psi &  & \sin\psi\cos\theta &  & \cos\psi\sin\theta
\end{array}\right)\left(\begin{array}{c}
\dot{\varphi}\\
\dot{\theta}\\
\dot{\psi}
\end{array}\right) & \,\, & \text{(see (\ref{eq:nonholonomous:constraints}))}\\
\left(\begin{array}{c}
\dot{\varphi}\\
\dot{\theta}\\
\dot{\psi}
\end{array}\right) & = & \left(\begin{array}{c}
\dot{\varphi}\\
\dot{\theta}\\
\dot{\psi}
\end{array}\right)\\
\left(\begin{array}{c}
\ddot{\varphi}\\
\ddot{\theta}\\
\ddot{\psi}
\end{array}\right) & = & \left(\begin{array}{c}
2\dot{\varphi}\dot{\theta}\tan\theta+\frac{5}{3}\dot{\theta}\dot{\psi}\cos\theta\\
\frac{4}{5r}g\sin\theta-\frac{6}{5}\dot{\varphi}\dot{\psi}\cos\theta+\frac{1}{2}\dot{\psi}^{2}\sin(2\theta)\\
\frac{2\dot{\varphi}\dot{\theta}}{\cos\theta}
\end{array}\right) &  & \text{(see (\ref{eq:lam:ddq}))}
\end{array}\label{eq:disk:state:eq}
\end{equation}

\section{Illustrations}

\begin{figure}[h]
\centering\includegraphics[width=9cm]{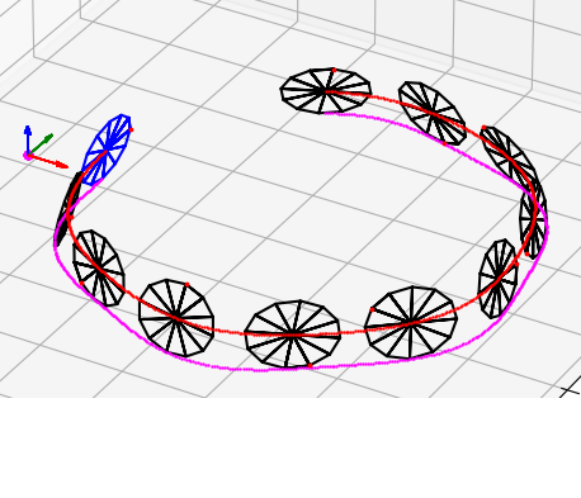}\caption{Disk rolling on a plane with precession}

\label{fig:rollingdisksimu1}
\end{figure}

Let us simulate the rolling disk with the following initial state
\begin{equation}
(c_{1},c_{2},\varphi,\theta,\psi,\dot{\varphi},\dot{\theta},\dot{\psi})=(2,0,0,0.1,0,2.5,0,0)
\end{equation}
for $t\in[0,10]$. The simulation generates Figure \ref{fig:rollingdisksimu1}
where the blue disk corresponds to the initial state. The behavior
is consistent with the intuition we could have for a rolling disk.
The precession effect is visible by the fact that the trajectory is
not perfectly circular.

To have a better understanding of the state model we have derived
for the rolling disk, let us consider different cases.

\textbf{Case 1}: The disk has a vertical motion. 

It means that $\theta=\dot{\theta}=\ddot{\theta}=0$. From (\ref{eq:disk:state:eq}),
we have
\begin{equation}
\left(\begin{array}{c}
\ddot{\varphi}\\
0\\
\ddot{\psi}
\end{array}\right)=\left(\begin{array}{c}
0\\
\dot{\varphi}\dot{\psi}\\
0
\end{array}\right).
\end{equation}
We have either 
\begin{itemize}
\item $\dot{\varphi}=0$: the contact point is static and the disk spins
around its vertical axis 
\item $\dot{\psi}=0$ : the disk moves straight forward
\end{itemize}
which is what we could have expected. 

\textbf{Case 2}: The disk is horizontal.

In this case, $\theta=\pm\frac{\pi}{2}$. The state equations are
not valid since we have a singularity. Indeed for this specific situation,
the heading and the spin is even not clearly defined.

\textbf{Case 3}. The disk has a circular trajectory. 

It means that $\ddot{\psi}=0$ with $\theta\neq0$ and $\dot{\varphi}\neq0$.
From (\ref{eq:disk:state:eq}), we get
\begin{equation}
\left(\begin{array}{c}
\ddot{\varphi}\\
0\\
0
\end{array}\right)=\left(\begin{array}{c}
2\dot{\varphi}\dot{\theta}\tan\theta+\frac{5}{3}\dot{\theta}\dot{\psi}\cos\theta\\
\frac{4}{5r}g\sin\theta-\frac{6}{5}\dot{\varphi}\dot{\psi}\cos\theta+\frac{1}{2}\dot{\psi}^{2}\sin(2\theta)\\
\dot{\theta}
\end{array}\right)
\end{equation}
\emph{i.e.},
\begin{equation}
\left(\begin{array}{c}
\ddot{\varphi}\\
0\\
0
\end{array}\right)=\left(\begin{array}{c}
0\\
\frac{4}{5r}g\sin\theta-\frac{6}{5}\dot{\varphi}\dot{\psi}\cos\theta+\frac{1}{2}\dot{\psi}^{2}\sin(2\theta)\\
\dot{\theta}
\end{array}\right).
\end{equation}

\begin{figure}[h]
\centering\includegraphics[width=9cm]{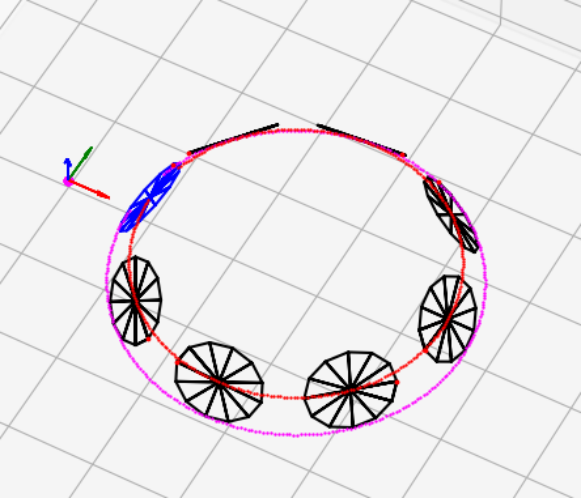}\caption{Disk rolling on a plane and performing a circle}

\label{fig:rollingdisksimu2}
\end{figure}

To illustrate this situation, let us simulate the rolling disk with
the following initial state
\begin{equation}
\left\{ \begin{array}{c}
(c_{1},c_{2},\varphi,\theta,\psi,\dot{\theta},\dot{\psi})=(2,0,0,0.5,0,0,1)\\
\dot{\varphi}=\frac{2}{3r\dot{\psi}}g\tan\theta+\frac{5}{6}\dot{\psi}\sin\theta
\end{array}\right.
\end{equation}
to satisfy the circular condition. For a time $t\in[0,6]$, the simulation
generates Figure \ref{fig:rollingdisksimu2} which corresponds indeed
to a circular motion.

\bigskip{}

The \texttt{Python} code associated to all examples can be found here: 
\begin{center}
\href{https://www.ensta-bretagne.fr/jaulin/rollingdisk.html}{https://www.ensta-bretagne.fr/jaulin/rollingdisk.html}
\par\end{center}

\bibliographystyle{plain}

\end{document}